\documentclass[11pt]{article}
\usepackage[left=1in,right=1in,top=1in,bottom=1in]{geometry}
\usepackage{times}
\usepackage{expl3}
\usepackage{cite}
\usepackage[table]{xcolor}
\usepackage{multirow}
\usepackage{stackengine} 
\usepackage{hhline}
\usepackage{lipsum}
\usepackage{titlesec}
\usepackage{breqn}
\usepackage{wrapfig}
\usepackage{epsfig}
\usepackage{graphics}
\usepackage{amsmath}
\usepackage[title]{appendix}
\usepackage{amssymb}
\usepackage{epstopdf}
\usepackage{boldline}
\usepackage{calligra}
\usepackage{url}
\usepackage{blindtext}
\usepackage{bm}

\newcommand{\define}{\stackrel{\mbox{\tiny def}}{=}}

\newtheorem{definition}{Definition}
\newtheorem{theorem}{Theorem}

\newtheorem{lemma}{Lemma}

\usepackage{mathtools}
\usepackage{epstopdf}
\usepackage{balance}
\usepackage{thmtools}
\usepackage{thm-restate}
\usepackage{hyperref}
\usepackage{cleveref}

\usepackage[ruled,vlined]{algorithm2e}
\include{pythonlisting}
\newcommand{\ostar}{\mathbin{\mathpalette\make@circled\star}}

\makeatletter
\newcommand{\removelatexerror}{\let\@latex@error\@gobble}
\makeatother
\makeatletter
\newcommand*{\rom}[1]{\expandafter\@slowromancap\romannumeral #1@}
\makeatother

\ExplSyntaxOn
\newcommand\latinabbrev[1]{
  \peek_meaning:NTF . {
    #1\@}%
  { \peek_catcode:NTF a {
      #1.\@ }%
    {#1.\@}}}
\ExplSyntaxOff

\titleclass{\subsubsubsection}{straight}[\subsubsection]

\begin{document}
\vspace{1cm}
\title{Generalized Hanson-Wright Inequality for Random Tensors}\vspace{1.8cm}
\author{Shih~Yu~Chang 
\thanks{Shih Yu Chang is with the Department of Applied Data Science,
San Jose State University, San Jose, CA, U. S. A. (e-mail: {\tt
shihyu.chang@sjsu.edu}).
           }}

\maketitle


\begin{abstract}
The Hanson-Wright inequality is an upper bound for tails of real quadratic forms in independent random variables. In this work, we extend the Hanson-Wright inequality for the Ky Fan $k$-norm for the polynomial function of the quadratic sum of random tensors under Einstein product. We decompose the quadratic tensors sum into the diagonal part and the coupling part. For the diagonal part, we can apply the generalized tensor Chernoff bound directly. But, for the coupling part,  we have to apply decoupling method first, i.e., decoupling inequality to bound expressions with dependent random tensors with independent random tensors, before applying generalized tensor Chernoff bound again to get the the tail probability of the Ky Fan $k$-norm of the coupling part sum of independent random tensors. At the end, the generalized Hanson-Wright inequality for the Ky Fan $k$-norm for the polynomial function of the quadratic sum of random tensors  can be obtained by the combination of the bound from the diagonal sum part and the bound from the coupling sum part.
\end{abstract}

\begin{keywords}
Hanson-Wright inequality, decoupling method, Ky Fan $k$-norm, Generalized Tensor Chernoff Bound, Hermitian tensors, Einstein product.
\end{keywords}

\section{Introduction}\label{sec:Introduction}

The Hanson-Wright inequality provides us an upper bound for tails of real quadratic forms
in independent subgaussian random variables. We define a random variable $Y$ is a $\alpha$-subgaussian if for every $\theta > 0$, we have $\mathrm{Pr}( | Y | \geq \theta ) \leq 2 \exp ( - \frac{ \theta^2}{2 \beta^2 })$. The Hanson-Wright inequality states that for any sequence of independent mean zero $\alpha$-subgaussian random variables $Y_1, \cdots,Y_n$, and any symmetric matrix $\mathbf{A}=(a_{i, j})_{i,j \leq n}$, we have
\begin{eqnarray}\label{eq:HW random variable form}
\mathrm{Pr}\left( \left\vert  \sum\limits_{i,j = 1}^n  a_{i,j}\left(Y_i Y_j -  \mathbb{E}(Y_i Y_j    ) \right)                     \right\vert  \geq \theta \right) \leq 2 \exp\left( -\frac{1}{C}\min\left\{ \frac{\theta^2}{\beta^4 \left\Vert \mathbf{A} \right\Vert_{\mbox{\tiny{HS}}}} , \frac{\theta}{\beta^2 \left\Vert \mathbf{A} \right\Vert_{\mbox{\tiny{OP}}}} \right\}  \right),
\end{eqnarray}
where $\left\Vert \mathbf{A} \right\Vert_{\mbox{\tiny{HS}}}$ is defined as $\left(  \sum\limits_{i,j = 1}^n  |a_{i,j}|^2   \right)^{1/2}$, and $\left\Vert \mathbf{A} \right\Vert_{\mbox{\tiny{OP}}}$ is defined as $\max\limits_{\left\Vert \mathbf{y}\right\Vert \leq 1} \left\Vert \mathbf{A}\mathbf{y}  \right\Vert_2$. The bound in Eq.~\eqref{eq:HW random variable form} was essentially proved in~\cite{adamczak2015note} in the symmetric case and in~\cite{hanson1971bound} in the zero mean case. The Hanson-Wright inequality has been applied to numerous applications in high-dimensional probability
and statistics, as well as in random matrix theory~\cite{vershynin2018high}. For example, the estimation of bound in Eq.~\eqref{eq:HW random variable form} is applied to the theory of compressed sensing with circulant type matrices~\cite{krahmer2014suprema}. In~\cite{adamczak2015note}, they applied Hanson-Wright inequality to study the concentration properties for sample covariance operators corresponding to Banach space-valued Gaussian random variables.

\emph{Tensor} was first introduced by William Ron Hamilton in 1846 and later became known to scientists through the publication of Levi-Civita’s book The Absolute Differential Calculus~\cite{levi1977absolute}. Because of its organized representation of data format and ability to reduce the complexity of multidimensional arrays, tensor has been gradually applied in various science and technology fields, such as physics~\cite{dahl2007tensor}, numerical computations~\cite{guan2019numerical}, unsupervised separation of unknown mixtures of speech signals in~\cite{wu2010robust, mirsamadi2016generalized}, multichannel signal filtering in~\cite{muti2007survey}, MIMO (multi-input multi-output) code-division in communication systems~\cite{de2008constrained, chen2011new}, passive sensing~\cite{kibangou2009blind, ribeiro2019low}, network signal processing in~\cite{shen2020topology} and image processing in~\cite{jiang2020framelet}. In~\cite{HW_T_SYChang_2021}, we first attempt to generalize Hanson-Wright inequality for the maximum eigenvalue of the quadratic sum of random Hermitian tensors under Einstein product. We first prove Weyl inequality for tensors under Einstein product and apply this fact to separate the quadratic form of random Hermitian tensors into diagonal sum and coupling (non-diagonal) sum parts. For the diagonal part, we can apply Bernstein inequality to bound the tail probability of the maximum eigenvalue of the sum of independent random Hermitian tensors directly. For coupling sum part, we have to apply decoupling method first, i.e., decoupling inequality to bound expressions with dependent random Hermitian tensors with independent random Hermitian tensors, before applying Bernstein inequality again to bound the tail probability of the maximum eigenvalue of the coupling sum of independent random Hermitian tensors. Finally, the Hanson-Wright inequality for the maximum eigenvalue of the quadratic sum of random Hermitian tensors under Einstein product can be obtained by the combination of the bound from the diagonal sum part and the bound from the coupling (non-diagonal) sum part.

In this work, we generalize our previous works from~\cite{HW_T_SYChang_2021} by considering Hanson-Wright inequality for the polynomial function of the quadratic sum of tensors. Our approach about the separation the quadratic sum into two parts: diagonal part and coupling part, is the same with the method adopted in~\cite{HW_T_SYChang_2021}, however, we have to apply the result from~\cite{chang2021general} to deal with the tail bounds analysis for the function of the random tensor sum. Our main theorem can be stated as follows. 
\begin{theorem}\label{thm:GHW Inequality}
We define a vector of random tensors $\overline{\mathcal{X}} \in \mathbb{C}^{(n \times I_1 \times \cdots \times I_M) \times (I_1 \times \cdots \times I_M)}$ as:
\begin{eqnarray}\label{eq:vec X def main}
\overline{\mathcal{X}} = \begin{bmatrix}
           \mathcal{X}_{1} \\
           \mathcal{X}_{2} \\
           \vdots \\
           \mathcal{X}_{n}
         \end{bmatrix},
\end{eqnarray}
where random Hermitian tensors $\mathcal{X}_{i} \in \mathbb{C}^{(I_1 \times \cdots \times I_M) \times (I_1 \times \cdots \times I_M)} $ are independent random positive definite tensors for $1 \leq i \leq n$. We also require another fixed tensor $\overline{\overline{\mathcal{A}}} \in  \mathbb{C}^{(n \times I_1 \times \cdots \times I_M) \times (n \times I_1 \times \cdots \times I_M)}$, which is defined as:
\begin{eqnarray}\label{eq:matrix A def main}
\overline{\overline{\mathcal{A}}} = \begin{bmatrix}
           \mathcal{A}_{1,1} &  \mathcal{A}_{1,2} & \cdots & \mathcal{A}_{1, n}  \\
           \mathcal{A}_{2,1} &  \mathcal{A}_{2,2} & \cdots & \mathcal{A}_{2, n}  \\
           \vdots & \vdots &   \vdots & \vdots \\
           \mathcal{A}_{n,1} &  \mathcal{A}_{n,2} & \cdots & \mathcal{A}_{n, n}  \\
         \end{bmatrix},
\end{eqnarray}
where $\mathcal{A}_{i,j} \in \mathbb{C}^{(I_1 \times \cdots \times I_M) \times (I_1 \times \cdots \times I_M)}$ are Hermitian tensors also. We define the tensor $\mathcal{D}_i$ associated to diagonal part of the tensor $\overline{\overline{\mathcal{A}}}$ as:
\begin{eqnarray}\label{eq:diagonal tensor sum main}
\mathcal{D}_i = \mathcal{X}_i \star_M \mathcal{A}_{i,i} \star_M  \mathcal{X}_i~~\mbox{for $i \in \{1,2,\ldots,n\}$};
\end{eqnarray}
and $\mathcal{C}_{\j}$ represents tensors associated to non-diagonal part of the tensor $\overline{\overline{\mathcal{A}}}$, which is defined as:
\begin{eqnarray}\label{eq:nondiagonal tensor sum main}
\mathcal{C}_{\j} =  \mathcal{X}_i \star_M \mathcal{A}_{i, j} \star_M  \mathcal{X}_j~~\mbox{for $i, j \in \{1,2,\ldots,n\}$ but $i \neq j$},
\end{eqnarray}
where the index $\j$ is determined by a pair of indices $i, j$. We will assume that tensors $\mathcal{D}_i$ and $\mathcal{C}_{\j}$ are positive definite tensors. For any $i$, we also assume that 
\begin{eqnarray}\label{eq:tensor commute assumption main}
\mathcal{X}_i \star_M \left( \sum\limits_{\ell=1, \neq i}^{n} \mathcal{A}_{i, \ell}\star_M  \mathcal{X}_{\ell}  \right)
=  \left( \sum\limits_{\ell=1, \neq i}^{n} \mathcal{A}_{i, \ell}\star_M  \mathcal{X}_{\ell}  \right)  \star_M \mathcal{X}_i. 
\end{eqnarray}
Besides, we also require 
\begin{eqnarray}\label{eq1:lma:coup part tail bound main}
\left( \exp\left(t    \sum\limits_{\ell=1, \neq i}^{n} \mathcal{A}_{i, \ell}\star_M  \mathcal{X}_{\ell}   \right)\right)^j  \geq \exp\left(t  \left(   \sum\limits_{\ell=1, \neq i}^{n} \mathcal{A}_{i, \ell}\star_M  \mathcal{X}_{\ell}  \right)^j \right)~~\mbox{almost surely for any $t > 0$};
\end{eqnarray}
and, given a positive real number $\mathrm{R}_d$, we have
\begin{eqnarray}
\lambda_{\max}\left( \mathcal{X}_{i} \star_M \mathcal{A}_{i, i} \star_M  \mathcal{X}_{i}   \right) \leq \mathrm{R}_d
\mbox{~~ almost surely for any $i \in \{1,2,\ldots, n\}$;}
\end{eqnarray}
and, given a positive real number $\mathrm{R}_c$, we also have
\begin{eqnarray}
\lambda_{\max}\left( \mathcal{A}_{i, \ell}\star_M  \mathcal{X}_{\ell}   \right) \leq \mathrm{R}_c
\mbox{~~ almost surely for any $i, l \in \{1,2,\ldots, n\}$ and $i \neq l$.}
\end{eqnarray}
We also assume that there is a Ky Fan $k$-norm bound for the exponent $j$ (positive integer) of random tensors $\mathcal{X}_i$, which is
\begin{eqnarray}\label{eq:Ky Fan norm bound for tensor assumption main}
\left\Vert \mathcal{X}^j_i \right\Vert_{(k)} \leq \mathrm{K}_{i,j,k},
\end{eqnarray}
where $i \in \{1,2,\ldots, n\}$ and $\mathrm{K}_{i,j,k} > 0$. Given a real polynomial $f(x) = a_0 + a_1 x + \ldots + a_m x^m$ and assume that $\Theta - \left\vert a_0 \right\vert k  = \sum\limits_{j=1}^m \theta_j$ with $\theta_j > 0$, we then have
\begin{eqnarray}
\mathrm{Pr}\left(  \left\Vert f \left( \overline{\mathcal{X}}^{\mathrm{T}} 
\overline{\overline{\mathcal{A}}} \overline{\mathcal{X}} \right) \right\Vert_{(k)} \geq \Theta   \right) \leq ~~~~~~~~~~~~~~~~~~~~~~~~~~~~~~~~~~~~~~~~~~~~~~~~~~~~~~~~~~~~~~~~~~~~~~~~~~~~~~~~~~~~~~~~~~~~~~~~~~~~~~~~ \nonumber \\
\sum\limits_{j=1}^m \left[ \mathrm{Pr}\left(   \left\Vert  \left(  \sum\limits_{i=1}^{n} \mathcal{D}_i  \right)^j \right\Vert_{(k)}    \geq \frac{1}{2^j} \left( \frac{\theta_j}{   \left\vert a_j \right\vert    } \right)  \right)  + \mathrm{Pr}\left(   \left\Vert  \left(    \sum\limits_{\j=1}^{n^2 - n} \mathcal{C}_{\j} \right)^j \right\Vert_{(k)}    \geq \frac{1}{2^j} \left( \frac{\theta_j}{   \left\vert a_j \right\vert    } \right)  \right)\right] ~~~~~~~~~~~~ \nonumber \\
\leq \sum\limits_{j=1}^m \inf\limits_{t > 0} e^{- \frac{ \theta_j t }{2^j \left\vert a_j \right\vert    } } \Bigg\{   \sum\limits_{i=1}^n \frac{k}{n} \Big[ 1 +\left( e^{n \mathrm{R}_{d} t} - 1 \right) \mathbb{E}\left(\sigma_1(\mathcal{D}_i)\right)  
+   C_{\mbox{\tiny Cher}} \left( e^{n\mathrm{R}_{d} t} - 1 \right) \Xi(\mathcal{D}_i) \Big] \Bigg\} 
~~~~~~~~~~~~~~~~~ \nonumber \\
+ \sum\limits_{j=1}^m   D_2 \Bigg\{ \sum\limits_{i=1}^{n}  \inf\limits_{t >0} e^{-   \frac{\theta_j t }{ 2^{j} n^{j-1}  \left\vert a_j \right\vert D_2  \mathrm{K}_{i,j,k}  }}  \Big\{ \sum\limits_{\ell=1, \neq i}^n  \frac{k}{n-1}\ \Big[ 1 +\left( e^{(n-1)\mathrm{R}_{c} t} - 1 \right) \mathbb{E}\left(\sigma_1\left( \mathcal{A}_{i, \ell}\star_M  \mathcal{X}_{\ell}     \right) \right) \nonumber \\
+   C_{\mbox{\tiny Cher}} \left( e^{(n-1)\mathrm{R}_{c} t} - 1 \right) \Xi( \mathcal{A}_{i, \ell}\star_M  \mathcal{X}_{\ell}   ) \Big]  \Big\} \Bigg\}.
\end{eqnarray}
\end{theorem}


The rest of this paper is organized as follows. In Section~\ref{sec:Preliminary Tensor Concepts and Genearlized Tensor Bounds}, we review tensors under Einstein product and tail bounds for the random tensors sum obtained by the marjoization approach which will be used in the later sections. In Section~\ref{sec:Function of Quadratic Form for Random Tensors and Diagonal Part}, we will first discuss Ky Fan $k$-norm tail probability formulation for the function of quadratic form with random tensors and develop the probability bound for the diagonal sum part. The coupling sum part will be discussed at Section~\ref{sec:Coupling Sum of Random Tensors} based on the decoupling inequality investigated in Section~\ref{sec:Decoupling Inequality for Tail Probability}. The proof of the main result of this work: generalized Hanson-Wright Inequality for random tensors, is given in Section~\ref{sec:Proof for Genearlized Hanson-Wright Inequality}. Finally, concluding discussions are given by Section~\ref{sec:Conclusion}.

\section{Preliminary Tensor Concepts and Genearlized Tensor Bounds}\label{sec:Preliminary Tensor Concepts and Genearlized Tensor Bounds} 

\subsection{Preliminary Tensor Concepts}\label{sec:Preliminary Tensor Concepts}

Throughout this work, scalars are represented by lower-case letters (e.g., $d$, $e$, $f$, $\ldots$), vectors by boldfaced lower-case letters (e.g., $\bm{d}$, $\bm{e}$, $\bm{f}$, $\ldots$), matrices by boldfaced capitalized letters (e.g., $\bm{D}$, $\bm{E}$, $\bm{F}$, $\ldots$), and tensors by calligraphic letters (e.g., $\mathcal{D}$, $\mathcal{E}$, $\mathcal{F}$, $\ldots$), respectively. Tensors are multiarrays of values which are higher-dimensional generalizations from vectors and matrices. Given a positive integer $N$, let $[N] \define \{1, 2, \cdots ,N\}$. An \emph{order-$N$ tensor} (or \emph{$N$-th order tensor}) denoted by $\mathcal{X} \define (x_{i_1, i_2, \cdots, i_N})$, where $1 \leq i_j = 1, 2, \ldots, I_j$ for $j \in [N]$, is a multidimensional array containing $\prod_{n=1}^N I_n$ entries. 
Let $\mathbb{C}^{I_1 \times \cdots \times I_N}$ and $\mathbb{R}^{I_1 \times \cdots \times I_N}$ be the sets of the order-$N$ $I_1 \times \cdots \times I_N$ tensors over the complex field $\mathbb{C}$ and the real field $\mathbb{R}$, respectively. For example, $\mathcal{X} \in \mathbb{C}^{I_1 \times \cdots \times I_N}$ is an order-$N$ multiarray, where the first, second, ..., and $N$-th dimensions have $I_1$, $I_2$, $\ldots$, and $I_N$ entries, respectively. Thus, each entry of $\mathcal{X}$ can be represented by $x_{i_1, \cdots, i_N}$. For example, when $N = 3$, $\mathcal{X} \in \mathbb{C}^{I_1 \times I_2 \times I_3}$ is a third-order tensor containing entries $x_{i_1, i_2, i_3}$'s.

Without loss of generality, one can partition the dimensions of a tensor into two groups, say $M$ and $N$ dimensions, separately. Thus, for two order-($M$+$N$) tensors: $\mathcal{X} \define (x_{i_1, \cdots, i_M, j_1, \cdots,j_N}) \in \mathbb{C}^{I_1 \times \cdots \times I_M\times
J_1 \times \cdots \times J_N}$ and $\mathcal{Y} \define (y_{i_1, \cdots, i_M, j_1, \cdots,j_N}) \in \mathbb{C}^{I_1 \times \cdots \times I_M\times
J_1 \times \cdots \times J_N}$, according to~\cite{MR3913666}, the \emph{tensor addition} $\mathcal{X} + \mathcal{Y}\in \mathbb{C}^{I_1 \times \cdots \times I_M\times
J_1 \times \cdots \times J_N}$ is given by 
\begin{eqnarray}\label{eq: tensor addition definition}
(\mathcal{X} + \mathcal{Y} )_{i_1, \cdots, i_M, j_1 , \cdots , j_N} &\define&
x_{i_1, \cdots, i_M, j_1 , \cdots , j_N} \nonumber \\
& &+ y_{i_1, \cdots, i_M, j_1 , \cdots , j_N}. 
\end{eqnarray}
On the other hand, for tensors $\mathcal{X} \define (x_{i_1, \cdots, i_M, j_1, \cdots,j_N}) \in \mathbb{C}^{I_1 \times \cdots \times I_M\times
J_1 \times \cdots \times J_N}$ and $\mathcal{Y} \define (y_{j_1, \cdots, j_N, k_1, \cdots,k_L}) \in \mathbb{C}^{J_1 \times \cdots \times J_N\times K_1 \times \cdots \times K_L}$, according to~\cite{MR3913666}, the \emph{Einstein product} (or simply referred to as \emph{tensor product} in this work) $\mathcal{X} \star_{N} \mathcal{Y} \in  \mathbb{C}^{I_1 \times \cdots \times I_M\times
K_1 \times \cdots \times K_L}$ is given by 
\begin{eqnarray}\label{eq: Einstein product definition}
\lefteqn{(\mathcal{X} \star_{N} \mathcal{Y} )_{i_1, \cdots, i_M,k_1 , \cdots , k_L} \define} \nonumber \\ &&\sum\limits_{j_1, \cdots, j_N} x_{i_1, \cdots, i_M, j_1, \cdots,j_N}y_{j_1, \cdots, j_N, k_1, \cdots,k_L}. 
\end{eqnarray}
Note that we will often abbreviate a tensor product $\mathcal{X} \star_{N} \mathcal{Y}$ to ``$\mathcal{X} \hspace{0.05cm}\mathcal{Y}$'' for notational simplicity in the rest of the paper. 
This tensor product will be reduced to the standard matrix multiplication as $L$ $=$ $M$ $=$ $N$ $=$ $1$. Other simplified situations can also be extended as tensor–vector product ($M >1$, $N=1$, and $L=0$) and tensor–matrix product ($M>1$ and $N=L=1$). In analogy to matrix analysis, we define some basic tensors and elementary tensor operations as follows. 

\begin{definition}\label{def: zero tensor}
A tensor whose entries are all zero is called a \emph{zero tensor}, denoted by $\mathcal{O}$. 
\end{definition}

\begin{definition}\label{def: identity tensor}
An \emph{identity tensor} $\mathcal{I} \in  \mathbb{C}^{I_1 \times \cdots \times I_N\times
J_1 \times \cdots \times J_N}$ is defined by 
\begin{eqnarray}\label{eq: identity tensor definition}
(\mathcal{I})_{i_1 \times \cdots \times i_N\times
j_1 \times \cdots \times j_N} \define \prod_{k = 1}^{N} \delta_{i_k, j_k},
\end{eqnarray}
where $\delta_{i_k, j_k} \define 1$ if $i_k  = j_k$; otherwise $\delta_{i_k, j_k} \define 0$.
\end{definition}

In order to define \emph{Hermitian} tensor, the \emph{conjugate transpose operation} (or \emph{Hermitian adjoint}) of a tensor is specified as follows.  
\begin{definition}\label{def: tensor conjugate transpose}
Given a tensor $\mathcal{X} \define (x_{i_1, \cdots, i_M, j_1, \cdots,j_N}) \in \mathbb{C}^{I_1 \times \cdots \times I_M\times J_1 \times \cdots \times J_N}$, its conjugate transpose, denoted by
$\mathcal{X}^{H}$, is defined by
\begin{eqnarray}\label{eq:tensor conjugate transpose definition}
(\mathcal{X}^H)_{ j_1, \cdots,j_N,i_1, \cdots, i_M}  \define  
 x^*_{i_1, \cdots, i_M,j_1, \cdots,j_N},
\end{eqnarray}
where the star $*$ symbol indicates the complex conjugate of the complex number $x_{i_1, \cdots, i_M,j_1, \cdots,j_N}$. If a tensor $\mathcal{X}$ satisfies $ \mathcal{X}^H = \mathcal{X}$, then $\mathcal{X}$ is a \emph{Hermitian tensor}. 
\end{definition}
We will use symbol $\iota$ to represent $\sqrt{-1}$.

Following definition is about untiary tensors.
\begin{definition}\label{def: unitary tensor}
Given a tensor $\mathcal{U} \define (u_{i_1, \cdots, i_M, i_1, \cdots,i_M}) \in \mathbb{C}^{I_1 \times \cdots \times I_M\times I_1 \times \cdots \times I_M}$, if
\begin{eqnarray}\label{eq:unitary tensor definition}
\mathcal{U}^H \star_M \mathcal{U} = \mathcal{U} \star_M \mathcal{U}^H = \mathcal{I} \in \mathbb{C}^{I_1 \times \cdots \times I_M\times I_1 \times \cdots \times I_M},
\end{eqnarray}
then $\mathcal{U}$ is a \emph{unitary tensor}. 
\end{definition}
In this work, the symbol $\mathcal{U}$ is reserved for a unitary tensor. 

\begin{definition}\label{def: inverse of a tensor}
Given a \emph{square tensor} $\mathcal{X} \define (x_{i_1, \cdots, i_M, j_1, \cdots,j_M}) \in \mathbb{C}^{I_1 \times \cdots \times I_M\times I_1 \times \cdots \times I_M}$, if there exists $\mathcal{Y} \in \mathbb{C}^{I_1 \times \cdots \times I_M\times I_1 \times \cdots \times I_M}$ such that 
\begin{eqnarray}\label{eq:tensor invertible definition}
\mathcal{Y} \star_M \mathcal{X} = \mathcal{X} \star_M \mathcal{Y} = \mathcal{I},
\end{eqnarray}
then $\mathcal{Y}$ is the \emph{inverse} of $\mathcal{X}$. We usually write $\mathcal{Y} \define \mathcal{X}^{-1}$ thereby. 
\end{definition}

We also list other crucial tensor operations here. The \emph{trace} of a square tensor is equivalent to the summation of all diagonal entries such that 
\begin{eqnarray}\label{eq: tensor trace def}
\mathrm{Tr}(\mathcal{X}) \define \sum\limits_{1 \leq i_j \leq I_j,\hspace{0.05cm}j \in [M]} \mathcal{X}_{i_1, \cdots, i_M,i_1, \cdots, i_M}.
\end{eqnarray}
The \emph{inner product} of two tensors $\mathcal{X}$, $\mathcal{Y} \in \mathbb{C}^{I_1 \times \cdots \times I_M\times J_1 \times \cdots \times J_N}$ is given by 
\begin{eqnarray}\label{eq: tensor inner product def}
\langle \mathcal{X}, \mathcal{Y} \rangle \define \mathrm{Tr}\left(\mathcal{X}^H \star_M \mathcal{Y}\right).
\end{eqnarray}
According to Eq.~\eqref{eq: tensor inner product def}, the \emph{Frobenius norm} of a tensor $\mathcal{X}$ is defined by 
\begin{eqnarray}\label{eq:Frobenius norm}
\left\Vert \mathcal{X} \right\Vert \define \sqrt{\langle \mathcal{X}, \mathcal{X} \rangle}.
\end{eqnarray}

Now, we wish to state a lemma about Ky Fan $k$-norm of two complex tensors. The Ky Fan $k$-norm of a complex tensor $\mathcal{A} \in \mathbb{C}^{I_1 \times \cdots \times I_M \times J_1 \times \cdots \times J_N}$, denoted as $\left\Vert \mathcal{A} \right\Vert_{(k)}$, is the summation of the largest $k$ singular values defined as:
\begin{eqnarray}\label{eq:Ky Fan k-norm def}
\left\Vert \mathcal{A} \right\Vert_{(k)} = \sum\limits_{i=1}^{k} \sigma_i(\mathcal{A}),
\end{eqnarray}
where $\sigma_i(\mathcal{A})$ is the $i$-th largest singular value of the tensor $\mathcal{A}$. The singular values of a complex tensor $\mathcal{A}$ are values at the diagonal entries in the diagonal tensor of the tensor $\mathcal{A}$ after singular value decomposition (SVD), see Theorem 3.2 in~\cite{liang2019further}. We apply symbol $\bm{\sigma}(\mathcal{A})$ to represent a vector with length $\min (I_1 \times \cdots \times I_M,  J_1 \times \cdots \times J_N)$ that is composed of all singular values of the tensor $\mathcal{A}$.

Given two real vectors $\bm{a}$ and $\bm{b}$ with length $m$, we say that the vector $\bm{a}$ weakly majorizes the vector $\bm{b}$, written as $\bm{b} \prec_{w} \bm{a}$, if we have
\begin{eqnarray}\label{eq:weak major def}
\sum\limits_{i=1}^{k}b_i^{\downarrow} \leq \sum\limits_{i=1}^{k}a_i^{\downarrow} \mbox{~~for $k=1,2,\ldots,m$,}
\end{eqnarray} 
where $a_i^{\downarrow}$ and $b_i^{\downarrow}$ are descending sorted elements of the vectors $\bm{a}$ and $\bm{b}$, i.e., $a_1^{\downarrow} \geq a_2^{\downarrow} \geq \ldots \geq a_k^{\downarrow}$ and $b_1^{\downarrow} \geq b_2^{\downarrow} \geq \ldots \geq b_k^{\downarrow}$.

\begin{lemma}\label{lma:Ky Fan k-norm inequality}
Given two tensors $\mathcal{A}, \mathcal{B} \in \mathbb{C}^{I_1 \times \cdots \times I_M \times J_1 \times \cdots \times J_N}$, we have 
\begin{eqnarray}\label{eq1:lma:Ky Fan k-norm inequalitty}
\bm{\sigma}(\mathcal{A} + \mathcal{B}) \prec_{w} \bm{\sigma}(\mathcal{A}) + \bm{\sigma}(\mathcal{B});
\end{eqnarray}
and
\begin{eqnarray}\label{eq2:lma:Ky Fan k-norm inequalitty}
\left\Vert  \mathcal{A} + \mathcal{B} \right\Vert_{(k)} \leq \left\Vert  \mathcal{A} \right\Vert_{(k)}
+ \left\Vert  \mathcal{B} \right\Vert_{(k)}.
\end{eqnarray}
\end{lemma}
\textbf{Proof:}
From Theorem 3.2 in~\cite{liang2019further}, we will have corresponding matrices $\bm{A}$ and $\bm{B}$ to tensors $\mathcal{A}$ and $\mathcal{B}$, respectively. Since sets of singular values of tensors $\mathcal{A}$ and $\mathcal{B}$ agree wih sets of singular values of matrices $\bm{A}$ and $\bm{B}$, we have Eq.~\eqref{eq1:lma:Ky Fan k-norm inequalitty} from Theorem G.1.d. in~\cite{MR2759813},  

Eq.~\eqref{eq2:lma:Ky Fan k-norm inequalitty} is true from the definition of Eq.~\eqref{eq1:lma:Ky Fan k-norm inequalitty} and Ky Fan $k$-norm definition provided by Eq.~\eqref{eq:Ky Fan k-norm def}.
$\hfill \Box$


From Theorem 5.2 in~\cite{ni2019hermitian}, every Hermitian tensor $\mathcal{H} \in  \mathbb{C}^{I_1 \times \cdots \times I_N \times I_1 \times \cdots \times I_N}$ has following decomposition
\begin{eqnarray}\label{eq:Hermitian Eigen Decom}
\mathcal{H} &=& \sum\limits_{i=1}^r \lambda_i \mathcal{U}_i \star_1 \mathcal{U}^{H}_i, \mbox{
~with~~$\langle \mathcal{U}_i, \mathcal{U}_i \rangle =1$ and $\langle \mathcal{U}_i, \mathcal{U}_j \rangle = 0$ for $i \neq j$,}
\end{eqnarray}
where $\lambda_i \in \mathbb{R}$, $\mathcal{U}_i  \in  \mathbb{C}^{I_1 \times \cdots \times I_N \times 1}$ and $\mathcal{U}_i  \in  \mathbb{C}^{1 \times I_1 \times \cdots \times I_N}$. The values $\lambda_i$ are named as \emph{Hermitian eigevalues}, and the value of $r$ will be $I_1 \times \cdots \times I_N$, if we count multiplicities of $\lambda_i$. 

Consider a function $f: \mathbb{R} \rightarrow \mathbb{R}$, we define a map on diagonal tensors by applying the function to each diagonal entry. We can extend $f$ to a function on a Hermitian tensor $\mathcal{H} \in \mathbb{C}^{I_1 \times \cdots \times I_N \times I_1 \times \cdots \times I_N}$ using the eigenvalue decomposition:
\begin{eqnarray}\label{eq:Hermitian tensor func}
f(\mathcal{H}) &\define& \mathcal{U} \star_N f(\mathcal{S}) \star_N \mathcal{U}^H,
\end{eqnarray}
where $\mathcal{H}$ can be expressed as $\mathcal{H} = \mathcal{U} \star_N \mathcal{S} \star_N \mathcal{U}^{H}$ via SVD from Theorem 3.2 in~\cite{liang2019further}. The spectral mapping theorem says that each eigenvalue of $f(\mathcal{H})$ is equal to $f(\lambda)$ for some eigenvalue of the tensor $\mathcal{H}$. From Eq.~\eqref{eq:Hermitian tensor func}, we have the following relationship: 
\begin{eqnarray}
f(x) \leq g(x)~~\mbox{for $x \in [a, b]$} \Rightarrow  f(\mathcal{H}) \leq g(\mathcal{H})~~\mbox{
when the eigenvalue of $\mathcal{H}$ within  $[a, b]$}, 
\end{eqnarray}
where $f(\mathcal{H}) \leq g(\mathcal{H})$ indicates that the tensor obtained by $g(\mathcal{H}) - f(\mathcal{H})$ is a nonnegative definite tensor, i.e., all eigenvalues of the tensor $g(\mathcal{H}) - f(\mathcal{H})$ are nonnegative.  


We will present a lemma about a tensor inequality in the sense of tensor definitess. 
\begin{lemma}\label{lma:convexity inequality}
Given two positive definite tensors $\mathcal{A}$ and $\mathcal{B}$ with $n \geq 1$, we have
\begin{eqnarray}\label{eq1:lma:convexity inequality}
\left\Vert  (\mathcal{A} + \mathcal{B}  )^n\right\Vert_{(k)}^{1/n} \leq \left\Vert  \mathcal{A}^n \right\Vert_{(k)}^{1/n} + \left\Vert  \mathcal{B}^n\right\Vert_{(k)}^{1/n}. 
\end{eqnarray}
\end{lemma}
\textbf{Proof:}
From unfilding operation provided by Sec. 2.2 in~\cite{liang2019further}, the given two positive definite tensors $\mathcal{A}$ and $\mathcal{B}$ will be transformed into  two positive definite matrices $\bm{A}$ and $\bm{B}$, this lemma is proved from 
Corollary 3.17 in~\cite{aujla2003weak} since Ky Fan $k$-norm is an unitary invariant norm. 
$\hfill \Box$

\subsection{Generalized Tail Bounds for Random Tensors Sum}\label{sec:Generalized Tail Bounds for Random Tensors Summation}

From our previous work~\cite{chang2021general}, we have following theorem about generalized tensor Chernoff bound. We restate this bound by the following theorem. 
\begin{theorem}[Generalized Tensor Chernoff Bound]\label{thm:Generalized Tensor Chernoff Bound intro}
Consider a sequence $\{ \mathcal{X}_j  \in \mathbb{C}^{I_1 \times \cdots \times I_N  \times I_1 \times \cdots \times I_N} \}$ of independent, random, Hermitian tensors. Let $g$ be a polynomial function with degree $n$ and nonnegative coeffecients $a_0, a_1, \cdots, a_n$ raised by power $s \geq 1$, i.e., $g(x) = \left(a_0 + a_1 x  +\cdots + a_n x^n \right)^s$ with $s \geq 1$. Suppose following condition is satisfied:
\begin{eqnarray}\label{eq:special cond Chernoff Bound intro}
g \left( \exp\left(t \sum\limits_{j=1}^{m} \mathcal{X}_j \right)\right)  \geq \exp\left(t g \left( \sum\limits_{j=1}^{m} \mathcal{X}_j   \right) \right)~~\mbox{almost surely},
\end{eqnarray}
where $t > 0$. Moreover, we require  
\begin{eqnarray}
\mathcal{X}_i \geq \mathcal{O} \mbox{~~and~~} \lambda_{\max}(\mathcal{X}_j) \leq \mathrm{R}
\mbox{~~ almost surely.}
\end{eqnarray}
Then we have the following inequality:
\begin{eqnarray}\label{eq1:thm:GeneralizedTensorChernoffBound intro}
\mathrm{Pr} \left( \left\Vert g\left( \sum\limits_{j=1}^{m} \mathcal{X}_j  \right)\right\Vert_{(k)}  \geq \theta \right)  \leq  (n+1)^{s-1} \inf\limits_{t > 0} e^{- \theta t } \cdot ~~~~~~~~~~~~~~~~~~~~~~~~~~~~~~~~~~~~~~~~~~~~~~~~~~~~~~~~ \nonumber \\
 \left\{ ka_0^s + \sum\limits_{l=1}^{n} \sum\limits_{j=1}^m \frac{k a_l^{ls}}{m} \left[ 1 +\left( e^{mlsRt} - 1 \right) \mathbb{E}\left[\sigma_1(\mathcal{X}_j)\right]  +  C_{\mbox{\tiny Cher}} \left( e^{mlsRt} - 1 \right) \Xi(\mathcal{X}_j) \right] \right\},
\end{eqnarray}
where $C_{\mbox{\tiny Cher}}$ is a constant. We use $\mathcal{X}^*$ to represent a tensor obtained by taking  complex conjugate of each entry of the tensor $\mathcal{X}$ Let  $x_{i,j}$ and $y_{i,j}$ are entries of matrices  obtained from unfolded random real tensors $\frac{\mathcal{X} + \mathcal{X}^*}{2} -  \mathbb{E}\left( \frac{ \mathcal{X} +\mathcal{X}^*}{2} \right)$ and $\frac{\mathcal{X} - \mathcal{X}^*}{2} -   \mathbb{E}\left( \frac{ \mathcal{X} - \mathcal{X}^*}{2} \right)$, respectively. The matrices from unfolded tensors are obtained by the method presented in Section 2.2~\cite{liang2019further}. For notation simplicity, the term $\Xi(\mathcal{X})$ is defined as 
\begin{eqnarray}\label{eq:abbr of sigma 1 of a random tensor}
\Xi(\mathcal{X})&\define& \left[ \max\limits_i \left( \sum\limits_j \mathbb{E} x^2_{i,j}\right)^{1/2} + 
\max\limits_j \left( \sum\limits_i \mathbb{E} x^2_{i,j}\right)^{1/2} + \left( \sum\limits_{i,j} \mathbb{E} x^4_{i,j}\right)^{1/4} \right.   +  \nonumber \\
&  & \left.
 \max\limits_i \left( \sum\limits_j \mathbb{E} y^2_{i,j}\right)^{1/2} +  \max\limits_j \left( \sum\limits_i \mathbb{E} y^2_{i,j}\right)^{1/2} + \left( \sum\limits_{i,j} \mathbb{E} y^4_{i,j}\right)^{1/4} \right]. 
\end{eqnarray}
\end{theorem}


\section{Function of Quadratic Form for Random Tensors and Diagonal Part}\label{sec:Function of Quadratic Form for Random Tensors and Diagonal Part}

In this section, we will first discuss Ky Fan $k$-norm tail probability formulation for the function of quadratic form with random tensors in Section~\ref{sec:Quadratic Form Partition}. The probability bound for the diagonal sum part will be presented by Section~\ref{sec:Diagonal Part of Random Tensors Summation}. The coupling sum part will be discussed at next Section~\ref{sec:Coupling Sum of Random Tensors}.

\subsection{Quadratic Form Partition}\label{sec:Quadratic Form Partition}

We define a vector of random tensors $\overline{\mathcal{X}} \in \mathbb{R}^{(n \times I_1 \times \cdots \times I_M) \times (I_1 \times \cdots \times I_M)}$ as:
\begin{eqnarray}\label{eq:vec X def}
\overline{\mathcal{X}} = \begin{bmatrix}
           \mathcal{X}_{1} \\
           \mathcal{X}_{2} \\
           \vdots \\
           \mathcal{X}_{n}
         \end{bmatrix},
\end{eqnarray}
where random Hermitian tensors $\mathcal{X}_{i}$ are independent random tensors. We also require another fixed tensor $\overline{\overline{\mathcal{A}}} \in  \mathbb{R}^{(n \times I_1 \times \cdots \times I_M) \times (n \times I_1 \times \cdots \times I_M)}$, which is defined as:
\begin{eqnarray}\label{eq:matrix A def}
\overline{\overline{\mathcal{A}}} = \begin{bmatrix}
           \mathcal{A}_{1,1} &  \mathcal{A}_{1,2} & \cdots & \mathcal{A}_{1, n}  \\
           \mathcal{A}_{2,1} &  \mathcal{A}_{2,2} & \cdots & \mathcal{A}_{2, n}  \\
           \vdots & \vdots &   \vdots & \vdots \\
           \mathcal{A}_{n,1} &  \mathcal{A}_{n,2} & \cdots & \mathcal{A}_{n, n}  \\
         \end{bmatrix},
\end{eqnarray}
where $\mathcal{A}_{i,j} \in \mathbb{R}^{(I_1 \times \cdots \times I_M) \times (I_1 \times \cdots \times I_M)}$ are Hermitian tensors also.

By independence, we can represent $\overline{\mathcal{X}}^{\mathrm{T}} 
\overline{\overline{\mathcal{A}}} \overline{\mathcal{X}}$ as 
\begin{eqnarray}\label{eq:quadratic form diag sum and coupling sum}
\overline{\mathcal{X}}^{\mathrm{T}} 
\overline{\overline{\mathcal{A}}} \overline{\mathcal{X}}  &=& \sum\limits_{i=1,j=1}^{n} \mathcal{X}_i \star_M \mathcal{A}_{i,j} \star_M  \mathcal{X}_j  \nonumber \\
&=& \sum\limits_{i=1}^{n}  \mathcal{X}_i \star_M \mathcal{A}_{i,i} \star_M  \mathcal{X}_i +  \sum\limits_{1 \leq i \neq j \leq n} \mathcal{X}_i \star_M \mathcal{A}_{i, j} \star_M  \mathcal{X}_j  \nonumber \\
& \define &  \sum\limits_{i=1}^{n} \mathcal{D}_i +  \sum\limits_{\j=1}^{n^2 - n} \mathcal{C}_{\j},
\end{eqnarray}
where $\mathcal{D}_i$ represents tensors related to diagonal part of the tensor $\overline{\overline{\mathcal{A}}}$, which is defined as:
\begin{eqnarray}\label{eq:diagonal tensor sum}
\mathcal{D}_i = \mathcal{X}_i \star_M \mathcal{A}_{i,i} \star_M  \mathcal{X}_i ;
\end{eqnarray}
and $\mathcal{C}_{\j}$ represents tensors related to non-diagonal part of the tensor $\overline{\overline{\mathcal{A}}}$, which is defined as:
\begin{eqnarray}\label{eq:nondiagonal tensor sum}
\mathcal{C}_{\j} =  \mathcal{X}_i \star_M \mathcal{A}_{i, j} \star_M  \mathcal{X}_j.
\end{eqnarray}
We will assume that tensors $\mathcal{D}_i$ and $\mathcal{C}_{\j}$ are positive definite tensors. We further assume that 
\begin{eqnarray}\label{eq:tensor commute assumption}
\mathcal{X}_i \star_M \left( \sum\limits_{\ell=1, \neq i}^{n} \mathcal{A}_{i, \ell}\star_M  \mathcal{X}_{\ell}  \right)
=  \left( \sum\limits_{\ell=1, \neq i}^{n} \mathcal{A}_{i, \ell}\star_M  \mathcal{X}^{(2)}_{\ell}  \right)  \star_M \mathcal{X}_i. 
\end{eqnarray}

Given a polynomial $f(x) = a_0 + a_1 x + a_2 x^2 + \ldots + a_m x^m$ with real coeffecients, we have
\begin{eqnarray}\label{eq:func of Quad decomp1}
\left\Vert f \left( \overline{\mathcal{X}}^{\mathrm{T}} 
\overline{\overline{\mathcal{A}}} \overline{\mathcal{X}}  \right)    \right\Vert_{(k)} &=&
\left\Vert a_0 \mathcal{I} + a_1 \left( \overline{\mathcal{X}}^{\mathrm{T}} 
\overline{\overline{\mathcal{A}}} \overline{\mathcal{X}}  \right)  + \ldots + a_m \left( \overline{\mathcal{X}}^{\mathrm{T}} 
\overline{\overline{\mathcal{A}}} \overline{\mathcal{X}}  \right)^m  \right\Vert_{(k)} \nonumber \\
&=& \left\Vert a_0 \mathcal{I} + a_1 \left(  \sum\limits_{i=1}^{n} \mathcal{D}_i +  \sum\limits_{\j=1}^{n^2 - n} \mathcal{C}_{\j}     \right)  \right. \nonumber \\
&  & \left. + \ldots + a_m \left( \sum\limits_{i=1}^{n} \mathcal{D}_i +  \sum\limits_{\j=1}^{n^2 - n} \mathcal{C}_{\j}    \right)^m  \right\Vert_{(k)} \nonumber \\
& \leq & \left\vert a_0 \right\vert  k + \left\vert a_1 \right\vert \left\Vert   \sum\limits_{i=1}^{n} \mathcal{D}_i +  \sum\limits_{\j=1}^{n^2 - n} \mathcal{C}_{\j}  \right\Vert_{(k)} + \left\vert a_2 \right\vert \left\Vert   \left( \sum\limits_{i=1}^{n} \mathcal{D}_i +  \sum\limits_{\j=1}^{n^2 - n} \mathcal{C}_{\j} \right)^2 \right\Vert_{(k)}  \nonumber \\
&  & + \ldots + \left\vert a_m \right\vert \left\Vert \left(  \sum\limits_{i=1}^{n} \mathcal{D}_i +  \sum\limits_{\j=1}^{n^2 - n} \mathcal{C}_{\j} \right)^m \right\Vert_{(k)}.
\end{eqnarray}
Then, we can have the following tail probability bound for $\left\Vert f \left( \overline{\mathcal{X}}^{\mathrm{T}} 
\overline{\overline{\mathcal{A}}} \overline{\mathcal{X}}  \right)    \right\Vert_{(k)} $ based on Eq.~\eqref{eq:func of Quad decomp1}. This bound can be expressed as
\begin{eqnarray}\label{eq:tail prob decomp1}
\mathrm{Pr}\left(  \left\Vert  f \left( \overline{\mathcal{X}}^{\mathrm{T}} 
\overline{\overline{\mathcal{A}}} \overline{\mathcal{X}}  \right) \right\Vert_{(k)}   \geq \Theta \right) &=& 
\mathrm{Pr}\left( \left\Vert  f \left( \overline{\mathcal{X}}^{\mathrm{T}} 
\overline{\overline{\mathcal{A}}} \overline{\mathcal{X}}  \right) \right\Vert_{(k)} -  \left\vert a_0 \right\vert k  \geq \Theta -  \left\vert a_0 \right\vert k \right)  \nonumber \\
& \leq & \sum\limits_{j=1}^{m} \mathrm{Pr}\left(   \left\vert a_j \right\vert \left\Vert \left(  \sum\limits_{i=1}^{n} \mathcal{D}_i +  \sum\limits_{\j=1}^{n^2 - n} \mathcal{C}_{\j} \right)^j \right\Vert_{(k)}    \geq \theta_j \right) 
\end{eqnarray} 
where we assume that $\Theta - \left\vert a_0 \right\vert k  = \sum\limits_{j=1}^m \theta_j$ with $\theta_j > 0$. For any $\left\vert a_j \right\vert \neq 0$~\footnote{If any $ \left\vert a_j \right\vert = 0$, the summation of Eq.~\eqref{eq:tail prob decomp1} will skip the terms with $ \left\vert a_j \right\vert = 0$. }, we have the following inequality:
\begin{eqnarray}\label{eq:tail prob decomp for each j}
 \mathrm{Pr}\left(   \left\vert a_j \right\vert \left\Vert  \left( \sum\limits_{i=1}^{n} \mathcal{D}_i +  \sum\limits_{\j=1}^{n^2 - n} \mathcal{C}_{\j} \right)^j \right\Vert_{(k)}    \geq \theta_j \right) &=&  \mathrm{Pr}\left(   \left\Vert \left(  \sum\limits_{i=1}^{n} \mathcal{D}_i +  \sum\limits_{\j=1}^{n^2 - n} \mathcal{C}_{\j} \right)^j \right\Vert_{(k)}    \geq \frac{\theta_j}{   \left\vert a_j \right\vert    } \right)  \nonumber \\
 &=&  \mathrm{Pr}\left(   \left\Vert  \left(  \sum\limits_{i=1}^{n} \mathcal{D}_i +  \sum\limits_{\j=1}^{n^2 - n} \mathcal{C}_{\j} \right)^j \right\Vert^{1/j}_{(k)}    \geq \left( \frac{\theta_j}{   \left\vert a_j \right\vert    } \right)^{1/j}  \right)  \nonumber \\
& \leq_1 & \mathrm{Pr}\left(   \left\Vert  \left(  \sum\limits_{i=1}^{n} \mathcal{D}_i  \right)^j \right\Vert^{1/j}_{(k)}    \geq \frac{1}{2} \left( \frac{\theta_j}{   \left\vert a_j \right\vert    } \right)^{1/j}  \right)  \nonumber \\
&  & + \mathrm{Pr}\left(   \left\Vert  \left(    \sum\limits_{\j=1}^{n^2 - n} \mathcal{C}_{\j} \right)^j \right\Vert^{1/j}_{(k)}    \geq \frac{1}{2} \left( \frac{\theta_j}{   \left\vert a_j \right\vert    } \right)^{1/j}  \right)  \nonumber \\
& =  & \mathrm{Pr}\left(   \left\Vert  \left(  \sum\limits_{i=1}^{n} \mathcal{D}_i  \right)^j \right\Vert_{(k)}    \geq \frac{1}{2^j} \left( \frac{\theta_j}{   \left\vert a_j \right\vert    } \right)  \right)  \nonumber \\
&  & + \mathrm{Pr}\left(   \left\Vert  \left(    \sum\limits_{\j=1}^{n^2 - n} \mathcal{C}_{\j} \right)^j \right\Vert_{(k)}    \geq \frac{1}{2^j} \left( \frac{\theta_j}{   \left\vert a_j \right\vert    } \right)  \right),
\end{eqnarray}
where we utilize Lemma~\ref{lma:convexity inequality} again in the inequality $\leq_1$ since we have
\begin{eqnarray} 
\left\Vert  \left(  \sum\limits_{i=1}^{n} \mathcal{D}_i +  \sum\limits_{\j=1}^{n^2 - n} \mathcal{C}_{\j} \right)^j \right\Vert^{1/j}_{(k)} \leq  \left\Vert  \left(  \sum\limits_{i=1}^{n} \mathcal{D}_i  \right)^j \right\Vert^{1/j}_{(k)}  +   \left\Vert  \left(    \sum\limits_{\j=1}^{n^2 - n} \mathcal{C}_{\j} \right)^j \right\Vert^{1/j}_{(k)}.
\end{eqnarray}

\subsection{Diagonal Part of Random Tensors Sum}\label{sec:Diagonal Part of Random Tensors Summation}

From Eqs.~\eqref{eq:func of Quad decomp1} and~\eqref{eq:tail prob decomp for each j}, random tensors involving diagonal part of the tensor $\overline{\overline{\mathcal{A}}}$ are $\left(  \sum\limits_{i=1}^{n} \mathcal{D}_i  \right)^j$ for $j=1,2,\ldots,m$, we can apply Theorem~\ref{thm:Generalized Tensor Chernoff Bound intro} to get the following lemma about the tail probability of the random tensor $\left(  \sum\limits_{i=1}^{n} \mathcal{D}_i  \right)^j$. 

\begin{lemma}[Bound for $\left(  \sum\limits_{i=1}^{n} \mathcal{D}_i  \right)^j$]\label{lma:diag part tail bound}
Consider a sequence $\{ \mathcal{D}_i  \in \mathbb{C}^{I_1 \times \cdots \times I_M  \times I_1 \times \cdots \times I_M} \}$ of independent, random, positive definite tensors. Let $g(x)$ be a polynomial function with degree $j$ as $g(x) = x^j$. Suppose following condition is satisfied:
\begin{eqnarray}\label{eq:special cond Chernoff Bound intro}
\left( \exp\left(t \sum\limits_{i=1}^{n} \mathcal{D}_i \right)\right)^j  \geq \exp\left(t  \left( \sum\limits_{i=1}^{n} \mathcal{D}_i   \right)^j \right)~~\mbox{almost surely},
\end{eqnarray}
where $t > 0$. Moreover, we require  
\begin{eqnarray}
\lambda_{\max}(\mathcal{D}_i) \leq \mathrm{R}_d
\mbox{~~ almost surely.}
\end{eqnarray}
Then we have the following inequality:
\begin{eqnarray}\label{eq1:thm:GeneralizedTensorChernoffBound intro}
\mathrm{Pr} \left( \left\Vert \left( \sum\limits_{i=1}^{n} \mathcal{D}_i  \right)^j \right\Vert_{(k)}  \geq \frac{ \theta_j}{2^j \left\vert a_j \right\vert    } \right)  &\leq&   \inf\limits_{t > 0} e^{- \frac{ \theta_j t }{2^j \left\vert a_j \right\vert    } } \Bigg\{   \sum\limits_{i=1}^n \frac{k}{n} \Big[ 1 +\left( e^{n\mathrm{R}_{d} t} - 1 \right) \mathbb{E}\left(\sigma_1(\mathcal{D}_i)\right)  \nonumber \\
&  &   +   C_{\mbox{\tiny Cher}} \left( e^{n \mathrm{R}_{d} t} - 1 \right) \Xi(\mathcal{D}_i) \Big] \Bigg\},
\end{eqnarray}
where $C_{\mbox{\tiny Cher}}$ is a constant and $\Xi(\mathcal{D}_i)$ is determined by Eq.~\eqref{eq:abbr of sigma 1 of a random tensor}.
\end{lemma}

\section{Coupling Sum of Random Tensors}\label{sec:Coupling Sum of Random Tensors}

The purpose of this section is to consider the tail probability bound for $\left\Vert \left( \sum\limits_{\j = 1}^{n^2 - n} \mathcal{C}_{\j} \right)^j  \right\Vert_{(k)}$. However, different from $\mathcal{D}_i$, the random tensors among $\mathcal{C}_{\j}$ are not indepedent. In Section~\ref{sec:Decoupling Inequality for Tail Probability}, we will present a decoupling inequality for function of random tensors. This decoupling inequality will help us to derive the tail bound for the random tensor summation $\left\Vert \left( \sum\limits_{\j = 1 }^{n^2 - n} \mathcal{C}_{\j } \right)^j  \right\Vert_{(k)}$ in Section~\ref{sec:Coupling Part of Random Tensors Summation}.

\subsection{Decoupling Inequality for Tail Probability}\label{sec:Decoupling Inequality for Tail Probability}


Before presenting the decoupling inequality for dependent random tensors, we have to prepare several lemmas first. 

\begin{lemma}\label{lma:lemma1}
Let $\mathcal{X}, \mathcal{Y}$ be two independent and identically distributed random Hermitian tensors with $\mathbb{E}(\mathcal{X}) = \mathbb{E}(\mathcal{Y}) = \mathcal{O}$. Then
\begin{eqnarray}
\mathrm{Pr} \left( \left\Vert  \mathcal{X} \right\Vert_{(k)} \geq \theta \right) \leq 
3 \mathrm{Pr} \left( \left\Vert  \mathcal{X} + \mathcal{Y} \right\Vert_{(k)}  \geq \frac{2\theta}{3} \right),
\end{eqnarray}
where $\theta > 0$.
\end{lemma}
\textbf{Proof:}
Let $\mathcal{Z}$ be another independent and identically distributed random Hermitian tensors compared to random Hermitian tensors $\mathcal{X}, \mathcal{Y}$ with $\mathbb{E}(\mathcal{Z})  = \mathcal{O}$. Then, we have
\begin{eqnarray}
\mathrm{Pr} \left( \left\Vert  \mathcal{X} \right\Vert_{(k)} \geq \theta \right)  
&=& \mathrm{Pr} \left( \left\Vert  ( \mathcal{X} +  \mathcal{Y}) + 
( \mathcal{X} +  \mathcal{Z}) - ( \mathcal{Y} +  \mathcal{Z})  \right\Vert_{(k)} \geq 2 \theta \right) \nonumber \\
&\leq& \mathrm{Pr} \left( \left\Vert  \mathcal{X} +  \mathcal{Y} \right\Vert_{(k)} \geq \frac{2 \theta}{3}  \right)
+   \mathrm{Pr} \left( \left\Vert   \mathcal{X} +  \mathcal{Z} \right\Vert_{(k)}   \geq \frac{2 \theta}{3}  \right) \nonumber \\
&  & + 
 \mathrm{Pr} \left( \left\Vert  \mathcal{Y} +  \mathcal{Z}\right\Vert_{(k)}  \geq \frac{2 \theta}{3}   \right)
\nonumber \\
&=& 3    \mathrm{Pr} \left( \left\Vert  \mathcal{Y} +  \mathcal{Z} \right\Vert_{(k)}   \geq \frac{2 \theta}{3}   \right)
\end{eqnarray}
$\hfill \Box$

\begin{lemma}\label{lma:prop1}
Let $\mathcal{X} \in \mathfrak{B}$, where $\mathfrak{B}$ is the Banach space with spectral norm,  be any zero mean random Hermitian tensor. Then for all non-random Hermitian tensor $\mathcal{A}$ same dimensions with $\mathcal{X}$ and $ \left\Vert \mathcal{A} \right\Vert_{(k)}> 0$, we have 
\begin{eqnarray}
\mathrm{Pr}\left( \left\Vert  \mathcal{A} + \mathcal{X}  \right\Vert_{(k)} \geq \left\Vert \mathcal{A}  \right\Vert_{(k)}  \right) \geq \frac{1}{4} \inf\limits_{f \in F}\frac{  (\mathbb{E}(\left\vert f(\mathcal{X})\right\vert))^2   }{\mathbb{E}(f^2(\mathcal{X}))}
\end{eqnarray}
where $F$ is the family of linear functionals on $\mathfrak{B}$.
\end{lemma}
\textbf{Proof:}
Note that if $x$ is a random variable with $\mathbb{E}x = 0$, then we have $\mathrm{Pr}(x \geq 0) \geq \frac{1}{4} \frac{   (  \mathbb{E}\left\vert x \right\vert )^2 }{\mathbb{E} (x^2)}$. From this fact, we have
\begin{eqnarray}
\mathrm{Pr}\left(f(\mathcal{X}) \geq 0 \right) \geq   \frac{1}{4} \frac{  (\mathbb{E}(\left\vert f(\mathcal{X})\right\vert))^2   }{\mathbb{E}(f^2(\mathcal{X}))}.
\end{eqnarray} 
If $f \in F$ (norming functional) is such that $f(\mathcal{A}) = \left\Vert  \mathcal{A} \right\Vert_{(k)}$ and $\left\Vert f \right\Vert_{\mathfrak{B}'}$, where $\left\Vert f \right\Vert_{\mathfrak{B}'} = 1$ is the function norm with respect to the dual space of $\mathfrak{B}$, denoted as $\mathfrak{B}'$, we have $\{  \left\Vert  \mathcal{A} + \mathcal{X} \right\Vert_{(k)}   \geq  \left\Vert  \mathcal{A} \right\Vert_{(k)}     \}$ contains $\{  f \left( \mathcal{A} + \mathcal{X} \right) \geq f \left( \mathcal{A} \right) \} = \{ f \left( \mathcal{X} \right) \geq 0 \}$. 
$\hfill \Box$

\begin{lemma}\label{lma:lemma2}
Let $\mathcal{A}_{i_1}, \mathcal{A}_{i_1,i_2}, \mathcal{A}_{i_1,i_2,i_3},\ldots, \mathcal{A}_{i_1,i_2,\ldots,i_m}, \mathcal{B}$ be non-random Hermitian tensors, and let $\{\beta_i \}$ be a sequence of independent and symmetric Bernoulli random variables, that is $\mathrm{Pr}(\beta_i = 1) = \mathrm{Pr}(\beta_i = -1) = \frac{1}{2}$. Then, we have 
\begin{eqnarray}
\mathrm{Pr}\left( \left\Vert \mathcal{B} + \sum\limits_{j=1}^{m} \sum\limits_{1 \leq i_1 \neq i_2 \neq \ldots \neq i_j \leq n}  \mathcal{A}_{i_1,i_2,\ldots,i_j}\beta_{i_1}\beta_{i_2}\ldots \beta_{i_j}    \right\Vert_{(k)} \geq \left\Vert  \mathcal{B} \right\Vert_{(k)} \right) \geq C_m
\end{eqnarray}
where $C_m$ is a constant depend on $\mathcal{A}_{i_1}, \mathcal{A}_{i_1,i_2}, \mathcal{A}_{i_1,i_2,i_3},\ldots, \mathcal{A}_{i_1,i_2,\ldots,i_m}$, but independent of $\mathcal{B}$.
\end{lemma}
\textbf{Proof:}
By setting $\mathcal{X} = \sum\limits_{j=1}^{m} \sum\limits_{1 \leq i_1 \neq i_2 \neq \ldots \neq i_j \leq n}  \mathcal{A}_{i_1,i_2,\ldots,i_j}\beta_{i_1}\beta_{i_2}\ldots \beta_{i_j}$ and $\mathcal{A} = \mathcal{B}$ in Lemma~\ref{lma:prop1}, this lemma is proved. 
$\hfill \Box$

We are ready to present the main Theorem in this section about the bounds on the tail probability by the decoupling inequality.  
\begin{theorem}\label{thm:decoupling}
Let $\{\mathcal{X}_i\}$ be a sequence of indepdent random tensors and $\{\mathcal{X}^{(j)}_i\}$, $j=1,2\ldots,m$, be $m$ indepedent copies of $\{\mathcal{X}_i\}$. Also let $f_{i_1, i_2, \ldots,i_m}$ be families of tensor-valued function of $m$ variables. Then, for all $n \geq m \geq 2$ and $\theta > 0$, there exists a contant $D_m$ dependeing on $m$ only so that 
\begin{eqnarray}\label{eq1:thm:decoupling}
\mathrm{Pr}\left(  \left\Vert  \sum\limits_{1 \leq i_1 \neq i_2 \neq \ldots \neq i_m \leq n}  f_{i_1, i_2,\ldots,i_m}\left(\mathcal{X}^{(1)}_{i_1},  \mathcal{X}^{(1)}_{i_2}, \ldots, \mathcal{X}^{(1)}_{i_m}  \right) \right\Vert_{(k)}   > \theta \right) ~~~~~~~~~~~~~~~~~~~~~~~~~~~  \nonumber \\
\leq D_m \mathrm{Pr}\left( D_m \left\Vert  \sum\limits_{1 \leq i_1 \neq i_2 \neq \ldots \neq i_m \leq n}  f_{i_1, i_2,\ldots,i_m}\left(\mathcal{X}^{(1)}_{i_1},  \mathcal{X}^{(2)}_{i_2}, \ldots, \mathcal{X}^{(m)}_{i_m}  \right) \right\Vert_{(k)}   > \theta \right).
\end{eqnarray}
\end{theorem}
\textbf{Proof:}

The proof shown here is obtained by applying the argument used in the proof of the
bound in the bivariate case followed by an inductive argument.

Let $\{\rho_i \}$ be a sequence of independent and symmetric Bernoulli random variables independent of random Hermitian tensors $\{\mathcal{X}^{(1)}_i\}, \{\mathcal{X}^{(2)}_i\}$. Let $(\mathcal{Z}^{(1)}, \mathcal{Z}^{(2)}) = (\mathcal{X}^{(1)}, \mathcal{X}^{(2)}) $ if $\rho_i = 1$, and  $(\mathcal{Z}^{(1)}, \mathcal{Z}^{(2)}) = (\mathcal{X}^{(2)}, \mathcal{X}^{(1)}) $ if $\rho_i = -1$. If $m=2$, we have
\begin{eqnarray}\label{eq2:thm:decoupling}
2^2 f_{i_1, i_2}\left(\mathcal{Z}_{i_1}^{(1)}, \mathcal{Z}_{i_2}^{(2)} \right) = (1+ \rho_{i_1})(1+ \rho_{i_2})
f_{i_1, i_2}\left(\mathcal{X}_{i_1}^{(1)}, \mathcal{X}_{i_2}^{(2)} \right) + (1+ \rho_{i_1})(1- \rho_{i_2})
f_{i_1, i_2}\left(\mathcal{X}_{i_1}^{(1)}, \mathcal{X}_{i_2}^{(1)} \right) \nonumber \\
+ (1 - \rho_{i_1})(1+ \rho_{i_2})
f_{i_1, i_2}\left(\mathcal{Z}_{i_1}^{(2)}, \mathcal{Z}_{i_2}^{(2)} \right) + (1 - \rho_{i_1})(1- \rho_{i_2})
f_{i_1, i_2}\left(\mathcal{Z}_{i_1}^{(2)}, \mathcal{Z}_{i_2}^{(1)} \right), 
\end{eqnarray}
where the sign $+$ is selected if the superscript of $\mathcal{X}_i$ agrees with that of $\mathcal{Z}_i$, and the sign $-$ is selected if the superscript of $\mathcal{X}_i$ disagrees with that of $\mathcal{Z}_i$. We set $\mathcal{S}_{n,2}$ as 
\begin{eqnarray}
 \mathcal{S}_{n,2} = \sum\limits_{1 \leq i_1 \neq i_2 \leq n} \Big[ (1+ \rho_{i_1})(1+ \rho_{i_2})
f_{i_1, i_2}\left(\mathcal{X}_{i_1}^{(1)}, \mathcal{X}_{i_2}^{(2)} \right) + (1+ \rho_{i_1})(1- \rho_{i_2})
f_{i_1, i_2}\left(\mathcal{X}_{i_1}^{(1)}, \mathcal{X}_{i_2}^{(1)} \right)  \nonumber \\
+ (1 - \rho_{i_1})(1+ \rho_{i_2})
f_{i_1, i_2}\left(\mathcal{Z}_{i_1}^{(2)}, \mathcal{Z}_{i_2}^{(2)} \right) + (1 - \rho_{i_1})(1- \rho_{i_2})
f_{i_1, i_2}\left(\mathcal{Z}_{i_1}^{(2)}, \mathcal{Z}_{i_2}^{(1)} \right) \Big].
\end{eqnarray}
If we define $\mathfrak{P}_2$ as a realization of $\mathcal{X}^{(1)}_i$ and $\mathcal{X}^{(2)}_i$ for $1 \leq i \leq n$, we have
\begin{eqnarray}
 \mathcal{S}_{n,2} = 2^2 \sum\limits_{1 \leq i_1 \neq i_2 \leq n} \mathbb{E}\left(
f_{i_1, i_2}\left(\mathcal{Z}_{i_1}^{(1)}, \mathcal{Z}_{i_2}^{(2)} \right) \vert \mathfrak{P}_2 \right). 
\end{eqnarray}

For $m > 2$ and any $1 \leq l_1, l_2, \ldots, l_m \leq 2$, we have
\begin{eqnarray}\label{eq:(4)}
2^m f_{i_1,\ldots,i_m}\left(\mathcal{Z}_{i_1}^{(l_1)}, \ldots, \mathcal{Z}_{i_m}^{(l_m)} \right)
= ~~~~~~~~~~~~~~~~~~~~~~~~~~~~~~~~~~~~~~~~~~~~~~~~~~~~~~~~~~~~~~~~~~~~ \nonumber \\
\sum\limits_{1 \leq j_1,\ldots, j_m \leq 2}(1 \pm_{(l_1, j_1)} \rho_{i_1})\ldots (1 \pm_{(l_m, j_m)} \rho_{i_m})
f_{i_1,\ldots,i_m}\left(\mathcal{X}_{i_1}^{(j_1)}, \ldots, \mathcal{X}_{i_m}^{(j_m)} \right),
\end{eqnarray}
where $ \pm_{(l_p, j_p)}$ is $+$ if $l_p = J_p$, and $ \pm_{(l_p, j_p)}$ is $-$ if $l_p \neq J_p$ for $p=1,2,\ldots,m$. Then the extension of $ \mathcal{S}_{n,2}$ becomes
\begin{eqnarray}
\mathcal{S}_{n,m} =  \sum\limits_{1 \leq i_1 \neq \ldots \neq  i_m \leq n}\sum\limits_{1 \leq j_1,\ldots, j_m \leq 2} f_{i_1,\ldots,i_m}\left(\mathcal{X}_{i_1}^{(j_1)}, \ldots, \mathcal{X}_{i_m}^{(j_m)} \right),
\end{eqnarray}
and we also can express $\mathcal{S}_{n,m}$ in terms of $\mathfrak{P}_2$  as
\begin{eqnarray}
\mathcal{S}_{n,m} = 2^m \sum\limits_{1 \leq i_1 \neq \ldots \neq  i_m \leq n} \mathbb{E}\left( f_{i_1,\ldots,i_m}\left(\mathcal{Z}_{i_1}^{(l_1)}, \ldots, \mathcal{Z}_{i_m}^{(l_m)} \right)  \vert \mathfrak{P}_2 \right).
\end{eqnarray}

From Lemma~\ref{lma:prop1}, we have 
\begin{eqnarray}\label{eq:(5)}
\mathrm{Pr}\left( \left\Vert \sum\limits_{1 \leq i_1 \neq \ldots \neq  i_m \leq n}   f_{i_1,\ldots,i_m}\left(\mathcal{X}_{i_1}^{(1)}, \ldots, \mathcal{X}_{i_m}^{(1)} \right)  \right\Vert_{(k)}     \geq \theta \right) \leq ~~~~~~~~~~~~~~~~~~~~~~~~~~~~~~~~~~~~~~~~~~~~~~~~~~~~~~~~~~~~ \nonumber \\
3 \mathrm{Pr}\left( 3 \left\Vert  \sum\limits_{1 \leq i_1 \neq \ldots \neq  i_m \leq n}  \left[f_{i_1,\ldots,i_m}\left(\mathcal{X}_{i_1}^{(1)}, \ldots, \mathcal{X}_{i_m}^{(1)} \right)  +  f_{i_1,\ldots,i_m}\left(\mathcal{X}_{i_1}^{(2)}, \ldots, \mathcal{X}_{i_m}^{(2)} \right)   \right]  \right\Vert_{(k)} \geq 2\theta \right) =  \nonumber \\
3 \mathrm{Pr}\Bigg( 3 \left\Vert \mathcal{S}_{n,m} +   \sum\limits_{1 \leq i_1 \neq \ldots \neq  i_m \leq n}  \Bigg[f_{i_1,\ldots,i_m}\left(\mathcal{X}_{i_1}^{(1)}, \ldots, \mathcal{X}_{i_m}^{(1)} \right) \right.
~~~~~~~~~~~~~~~~~~~~~~~~~~~~~~~~~~~~~~~~~~~~~~~~~~~~~~~~~~~ \nonumber \\ \left. +  f_{i_1,\ldots,i_m}\left( \mathcal{X}_{i_1}^{(2)}, \ldots, \mathcal{X}_{i_m}^{(2)} \right)   \Bigg]  - \mathcal{S}_{n,m}\right\Vert_{(k)} \geq 2\theta \Bigg) \leq  ~~~~~~~~~~~~~~~~~~~~~~~~~~~~~~~~~~~~~~~~~~~~~~~~~~~~~~~~~ \nonumber \\~~~~~~
3 \mathrm{Pr}\left( 3 \left\Vert \mathcal{S}_{n,m}  \right\Vert_{(k)} \geq \theta \right) +
~~~~~~~~~~~~~~~~~~~~~~~~~~~~~~~~~~~~~~~~~~~~~~~~~~~~~~~~~~~~~~~~~~~~~~~~~~~~~~~~~~~~~~~~~~~~~~~~~~~~~~~~~~~~~~ \nonumber \\
3 \mathrm{Pr}\left( 3 \left\Vert  \sum\limits_{1 \leq i_1 \neq \ldots \neq  i_m \leq n}  \sum\limits_{\substack{1 \leq j_1, \ldots j_m \leq 2  \\  \mbox{\tiny remove same $j'$}}} f_{i_1,\ldots,i_m}\left( \mathcal{X}_{i_1}^{(j_1)}, \ldots, \mathcal{X}_{i_m}^{(j_m)} \right) \right\Vert_{(k)} \geq \theta \right)  \leq_1  ~~~~~~~~~~~~~~~~~~~~~~~~~ \nonumber \\
3 \mathrm{Pr}\left( 3 \left\Vert \mathcal{S}_{n,m}  \right\Vert_{(k)} \geq \theta \right) +
~~~~~~~~~~~~~~~~~~~~~~~~~~~~~~~~~~~~~~~~~~~~~~~~~~~~~~~~~~~~~~~~~~~~~~~~~~~~~~~~~~~~~~~~~~~~~~~~~~~~~~~~~~~~~~\nonumber \\
 \sum\limits_{\substack{1 \leq j_1, \ldots j_m \leq 2  \\  \mbox{\tiny remove same $j'$}}} E_m \mathrm{Pr}\left( E_m \left\Vert  \sum\limits_{1 \leq i_1 \neq \ldots \neq  i_m \leq n}  f_{i_1,\ldots,i_m}\left( \mathcal{X}_{i_1}^{(j_1)}, \ldots, \mathcal{X}_{i_m}^{(j_m)} \right) \right\Vert_{(k)} \geq \theta \right),~~~~~~~~~~~~~~~~~~~~~~
\end{eqnarray}
where the last inequality $\leq_1$ is obtained by the triangle inequality of Ky Fan norm and the tail bound property, and the constant $E_m$ is depend on the value $m$ only.  

Given any fixed $1 \leq l_1, \ldots, l_m \leq 2$ such that not all $l'$ are equal, from Lemma~\ref{lma:lemma2} and Eq.~\eqref{eq:(4)}, we have 
\begin{eqnarray}\label{eq:(6)}
\mathrm{Pr}\left(2^m \left\Vert    \sum\limits_{1 \leq i_1 \neq \ldots \neq  i_m \leq n} f_{i_1,\ldots,i_m}\left( \mathcal{Z}_{i_1}^{(l_1)}, \ldots, \mathcal{Z}_{i_m}^{(l_m)} \right) \right\Vert_{(k)} \geq \left\Vert \mathcal{S}_{n,m} \right\Vert_{(k)}  \vert   \mathfrak{P}_2   \right) \geq  C_m.
\end{eqnarray}
By integrating over $\{\left\Vert \mathcal{S}_{n,m} \right\Vert_{(k)} \geq \theta\}$ and apply the fact that $\left\{\left(\mathcal{X}_i^{(1)}, \mathcal{X}_i^{(2)}\right) \mbox{for $i=1,\ldots,n$} \right\}$ has the same distribution as 
$\left\{\left(\mathcal{Z}_i^{(1)}, \mathcal{Z}_i^{(2)}\right) \mbox{for $i=1,\ldots,n$} \right\}$, we have
\begin{eqnarray}\label{eq:(7)}
\mathrm{Pr}\left(2^m \left\Vert    \sum\limits_{1 \leq i_1 \neq \ldots \neq  i_m \leq n} f_{i_1,\ldots,i_m}\left( \mathcal{X}_{i_1}^{(l_1)}, \ldots, \mathcal{X}_{i_m}^{(l_m)} \right) \right\Vert_{(k)} \geq \theta \right) = ~~~~~~~~~~~~~~~~~~~~~~~~~~~~~~~~~~~~~~~~~~ \nonumber \\
\mathrm{Pr}\left(2^m \left\Vert    \sum\limits_{1 \leq i_1 \neq \ldots \neq  i_m \leq n} f_{i_1,\ldots,i_m}\left( \mathcal{Z}_{i_1}^{(l_1)}, \ldots, \mathcal{Z}_{i_m}^{(l_m)} \right) \right\Vert_{(k)} \geq \theta \right) \geq C_m \mathrm{Pr}\left( \left\Vert \mathcal{S}_{n,m} \right\Vert_{(k)} \geq \theta \right).
\end{eqnarray}

We assume that the decoupling inequality is valid for $1,2,\ldots,m-1$. From Eqs.~\eqref{eq:(5)} and~\eqref{eq:(7)}, then we have
\begin{eqnarray}
\mathrm{Pr}\left( \left\Vert \sum\limits_{1 \leq i_1 \neq \ldots \neq  i_m \leq n}   f_{i_1,\ldots,i_m}\left(\mathcal{X}_{i_1}^{(1)}, \ldots, \mathcal{X}_{i_m}^{(1)} \right)  \right\Vert_{(k)}     \geq \theta \right) \leq 3 \mathrm{Pr}\left( 3 \left\Vert \mathcal{S}_{n,m}  \right\Vert_{(k)} \geq \theta \right) \nonumber \\
+ 
 \sum\limits_{\substack{1 \leq j_1, \ldots j_m \leq 2  \\  \mbox{\tiny remove same $j'$}}} E_m \mathrm{Pr}\left( E_m \left\Vert  \sum\limits_{1 \leq i_1 \neq \ldots \neq  i_m \leq n}  f_{i_1,\ldots,i_m}\left( \mathcal{X}_{i_1}^{(j_1)}, \ldots, \mathcal{X}_{i_m}^{(j_m)} \right) \right\Vert_{(k)} \geq \theta \right) \nonumber \\
\leq \frac{3}{C_m} \mathrm{Pr}\left(3 \cdot 2^m \left\Vert    \sum\limits_{1 \leq i_1 \neq \ldots \neq  i_m \leq n} f_{i_1,\ldots,i_m}\left( \mathcal{X}_{i_1}^{(l_1)}, \ldots, \mathcal{X}_{i_m}^{(l_m)} \right) \right\Vert_{(k)} \geq \theta \right) ~~~~~~~~~~~~~~ \nonumber \\
+ 
\sum\limits_{\substack{1 \leq j_1, \ldots j_m \leq 2  \\  \mbox{\tiny remove same $j'$}}} E_m \mathrm{Pr}\left( E_m \left\Vert  \sum\limits_{1 \leq i_1 \neq \ldots \neq  i_m \leq n}  f_{i_1,\ldots,i_m}\left( \mathcal{X}_{i_1}^{(j_1)}, \ldots, \mathcal{X}_{i_m}^{(j_m)} \right) \right\Vert_{(k)} \geq \theta \right) \nonumber \\
\leq_1 \sum\limits_{\substack{1 \leq j_1, \ldots j_m \leq 2  \\  \mbox{\tiny remove same $j'$}}} \overline{E}_m \mathrm{Pr}\left( \overline{E}_m \left\Vert  \sum\limits_{1 \leq i_1 \neq \ldots \neq  i_m \leq n}  f_{i_1,\ldots,i_m}\left( \mathcal{X}_{i_1}^{(j_1)}, \ldots, \mathcal{X}_{i_m}^{(j_m)} \right) \right\Vert_{(k)} \geq \theta \right) \nonumber \\
\leq_2  D_m \mathrm{Pr}\left( D_m \left\Vert  \sum\limits_{1 \leq i_1 \neq i_2 \neq \ldots \neq i_m \leq n}  f_{i_1, i_2,\ldots,i_m}\left(\mathcal{X}^{(1)}_{1},  \mathcal{X}^{(2)}_{2}, \ldots, \mathcal{X}^{(m)}_{m}  \right) \right\Vert_{(k)}   > \theta \right),~
\end{eqnarray}
where the inequality $\leq_1$ is obtained by adjusting contants $E_m$, and the inequality $\leq_2$ is obtained by  decoupling result for U-statistics of orders $2,\ldots, m-1$ of the induction. This theorem is proved.
$\hfill \Box$

Above proof method is extended from~\cite{de1995decoupling} to tensors, but we try to show those different bounding constants, lile $E_m, \overline{E}_m, D_m$, which are treated as same symbols in the original proof argument. This is misleading. 

\subsection{Coupling Part of Random Tensors Summation}\label{sec:Coupling Part of Random Tensors Summation}


From Eqs.~\eqref{eq:func of Quad decomp1} and~\eqref{eq:tail prob decomp for each j}, random tensors involving coupling part of the tensor $\overline{\overline{\mathcal{A}}}$ are $\left(  \sum\limits_{\j=1}^{n^2 - n} \mathcal{C}_{\j}  \right)^j$ for $j=1,2,\ldots,m$, we can apply Theorem~\ref{thm:Generalized Tensor Chernoff Bound intro} and Theorem~\ref{thm:decoupling} to get the following lemma about the tail probability of the random tensor $\left(  \sum\limits_{\j=1}^{n^2 - n} \mathcal{C}_{\j}  \right)^j$. 

\begin{lemma}[Bound for $\left(  \sum\limits_{\j=1}^{n^2 - n} \mathcal{C}_{\j}  \right)^j$]\label{lma:coup part tail bound}
Consider a sequence $\{ \mathcal{C}_{\j}  \in \mathbb{C}^{I_1 \times \cdots \times I_M  \times I_1 \times \cdots \times I_M} \}$ of random, positive definite tensors. Let $g(x)$ be a polynomial function with degree $j$ as $g(x) = x^j$, and let the random tensor $\tilde{\mathcal{C}}_{\j}$ be transformed from the random tensor $\mathcal{C}_{\j}$
as:
\begin{eqnarray}
\tilde{\mathcal{C}}_{\j} =  \mathcal{X}^{(1)}_i \star_M \mathcal{A}_{i, j} \star_M  \mathcal{X}^{(2)}_j  -  \mathbb{E}\left(\mathcal{X}^{(1)}_i\right) \star_M \mathcal{A}_{i,i} \star_M  \mathbb{E}\left( \mathcal{X}^{(2)}_j \right),
\end{eqnarray}
where the random tensors $\mathcal{X}^{(1)}_i, \mathcal{X}^{(2)}_j$ are copies from the random tensors $\mathcal{X}_i, \mathcal{X}_j$. Suppose the following condition is satisfied for any $i \in \{1,2,\ldots, n \}$:
\begin{eqnarray}\label{eq1:lma:coup part tail bound}
\left( \exp\left(t    \sum\limits_{\ell=1, \neq i}^{n} \mathcal{A}_{i, \ell}\star_M  \mathcal{X}_{\ell}   \right)\right)^j  \geq \exp\left(t  \left(   \sum\limits_{\ell=1, \neq i}^{n} \mathcal{A}_{i, \ell}\star_M  \mathcal{X}_{\ell}  \right)^j \right)~~\mbox{almost surely},
\end{eqnarray}
where $t > 0$. Moreover, we require  
\begin{eqnarray}
\lambda_{\max}\left( \mathcal{A}_{i, \ell}\star_M  \mathcal{X}_{\ell}   \right) \leq \mathrm{R}_c
\mbox{~~ almost surely for any $i, l \in \{1,2,\ldots, n\}$;}
\end{eqnarray}
and a Ky Fan bound for random tensor exponent, which is
\begin{eqnarray}\label{eq:Ky Fan norm bound for tensor assumption}
\left\Vert \mathcal{X}^j_i \right\Vert_{(k)} \leq \mathrm{K}_{i,j,k},
\end{eqnarray}
where $\mathrm{K}_{i,j,k} > 0$.

Then we have the following inequality:
\begin{eqnarray}\label{eq2:lma:coup part tail bound}
\mathrm{Pr}\left(   \left\Vert  \left(    \sum\limits_{\j=1}^{n^2 - n} \mathcal{C}_{\j} \right)^j \right\Vert_{(k)}    \geq \frac{1}{2^j} \left( \frac{\theta_j}{   \left\vert a_j \right\vert    } \right)  \right)
\leq D_2 \sum\limits_{i=1}^{n}  \inf\limits_{t >0} e^{-   \frac{\theta_j t }{ 2^{j} n^{j-1}  \left\vert a_j \right\vert D_2  \mathrm{K}_{i,j,k}  }} \nonumber \\
\times \left\{ \sum\limits_{\ell=1, \neq i}^n  \frac{k}{n-1}\ \Big[ 1 +\left( e^{(n-1)\mathrm{R}_{c} t} - 1 \right) \mathbb{E}\left(\sigma_1\left( \mathcal{A}_{i, \ell}\star_M  \mathcal{X}_{\ell}     \right) \right) \right. \nonumber \\
  \left. +   C_{\mbox{\tiny Cher}} \left( e^{(n-1)\mathrm{R}_{c} t} - 1 \right) \Xi( \mathcal{A}_{i, \ell}\star_M  \mathcal{X}_{\ell}   ) \Big]  \right\}.,
\end{eqnarray}
where $C_{\mbox{\tiny Cher}}$ is a constant and $\Xi( \mathcal{A}_{i, \ell}\star_M  \mathcal{X}^{(2)}_{\ell}   )$ is determined by Eq.~\eqref{eq:abbr of sigma 1 of a random tensor}, and $D_2$ comes from Theorem~\ref{thm:decoupling}.
\end{lemma}
\textbf{Proof:}
From Theorem~\ref{thm:decoupling}, we have 
\begin{eqnarray}\label{eq3:lma:coup part tail bound}
\mathrm{Pr}\left(   \left\Vert  \left(    \sum\limits_{\j=1}^{n^2 - n} \mathcal{C}_{\j} \right)^j \right\Vert_{(k)}    \geq \frac{1}{2^j} \left( \frac{\theta_j}{   \left\vert a_j \right\vert    } \right)  \right)
\leq D_2 \mathrm{Pr}\left( D_2   \left\Vert  \left(    \sum\limits_{\j=1}^{n^2 - n} \tilde{\mathcal{C}}_{\j} \right)^j \right\Vert_{(k)}    \geq \frac{1}{2^j} \left( \frac{\theta_j}{   \left\vert a_j \right\vert    } \right)  \right) \nonumber \\
= D_2 \mathrm{Pr}\left(  \left\Vert  \left(    \sum\limits_{\j=1}^{n^2 - n} \tilde{\mathcal{C}}_{\j} \right)^j \right\Vert_{(k)}    \geq \frac{\theta_j}{ 2^j  \left\vert a_j \right\vert D_2   }   \right) ~~~~~~~~~~~~~~~~~~~~~~~~~~~~~~~~~~~~~~~~~~~~~~~~~~~~~~~~~  \nonumber \\
 \leq_1  D_2 \sum\limits_{i=1}^{n} \mathrm{Pr}\left(  \left\Vert \left[ \mathcal{X}^{(1)}_i \star_M \left( \sum\limits_{\ell=1, \neq i}^{n} \mathcal{A}_{i, \ell}\star_M  \mathcal{X}^{(2)}_{\ell}  \right)    \right]^j     \right\Vert_{(k)}   \geq   \frac{\theta_j}{ 2^{j} n^{j-1}  \left\vert a_j \right\vert D_2   }    \right)
\end{eqnarray}
where we apply Lemma 4.3 from \cite{chang2021tensor} in $\leq_1$. From assumptions provided by Eq.~\eqref{eq:tensor commute assumption}, Eq~\eqref{eq:Ky Fan norm bound for tensor assumption} and the fact that $\left\Vert \mathcal{A} \star \mathcal{B}  \right\Vert_{(k)} \leq \left\Vert \mathcal{A}  \right\Vert_{(k)} \left\Vert \mathcal{B}  \right\Vert_{(k)} $, we can further bound each summand in Eq.~\eqref{eq3:lma:coup part tail bound} as 
\begin{eqnarray}\label{eq4:lma:coup part tail bound}
\mathrm{Pr}\left(  \left\Vert \left[ \mathcal{X}^{(1)}_i \star_M \left( \sum\limits_{\ell=1, \neq i}^{n} \mathcal{A}_{i, \ell}\star_M  \mathcal{X}^{(2)}_{\ell}  \right)    \right]^j     \right\Vert_{(k)}   \geq   \frac{\theta_j}{ 2^{j} n^{j-1}  \left\vert a_j \right\vert D_2   }    \right) ~~~~~~~~~~~~~~~~~~~~~~~~~~~~~~~~  \nonumber \\
= \mathrm{Pr}\left(  \left\Vert \left[ \left(\mathcal{X}^{(1)}_i \right)^j     \star_M \left( \sum\limits_{\ell=1, \neq i}^{n} \mathcal{A}_{i, \ell}\star_M  \mathcal{X}^{(2)}_{\ell}  \right)^j        \right] \right\Vert_{(k)}   \geq   \frac{\theta_j}{ 2^{j} n^{j-1}  \left\vert a_j \right\vert D_2   }    \right)  ~~~~~~~~~~~~~~~~~~~~~ \nonumber \\
\leq \mathrm{Pr}\left(  \left\vert  \left( \sum\limits_{\ell=1, \neq i}^{n} \mathcal{A}_{i, \ell}\star_M  \mathcal{X}^{(2)}_{\ell}  \right)^j  \right\Vert_{(k)}   \geq \frac{\theta_j}{ 2^{j} n^{j-1}  \left\vert a_j \right\vert D_2  \mathrm{K}_{i,j,k}  }  \right)  ~~~~~~~~~~~~~~~~~~~~~~~~~~~~~~~~~~~~~~~  \nonumber \\
\leq_1 \inf\limits_{t >0} e^{-   \frac{\theta_j t }{ 2^{j} n^{j-1}  \left\vert a_j \right\vert D_2  \mathrm{K}_{i,j,k}  }}
\left\{ \sum\limits_{\ell=1, \neq i}^n  \frac{k}{n-1}\ \Big[ 1 +\left( e^{(n-1) \mathrm{R}_{c} t} - 1 \right) \mathbb{E}\left(\sigma_1\left( \mathcal{A}_{i, \ell}\star_M  \mathcal{X}^{(2)}_{\ell}     \right) \right) \right. \nonumber \\
  \left. +   C_{\mbox{\tiny Cher}} \left( e^{(n-1)\mathrm{R}_{c}  t} - 1 \right) \Xi( \mathcal{A}_{i, \ell}\star_M  \mathcal{X}^{(2)}_{\ell}   ) \Big]  \right\}
\end{eqnarray}
where $\leq_1$ comes from Theorem~\ref{thm:Generalized Tensor Chernoff Bound intro}. 

Finally, from Eqs.~\eqref{eq3:lma:coup part tail bound} and~\eqref{eq4:lma:coup part tail bound}, we have
\begin{eqnarray}
\mathrm{Pr}\left(   \left\Vert  \left(    \sum\limits_{\j=1}^{n^2 - n} \mathcal{C}_{\j} \right)^j \right\Vert_{(k)}    \geq \frac{1}{2^j} \left( \frac{\theta_j}{   \left\vert a_j \right\vert    } \right)  \right)
\leq  D_2 \sum\limits_{i=1}^{n}  \inf\limits_{t >0} e^{-   \frac{\theta_j t }{ 2^{j} n^{j-1}  \left\vert a_j \right\vert D_2  \mathrm{K}_{i,j,k}  }} ~~~~~~ \nonumber \\
\times \left\{ \sum\limits_{\ell=1, \neq i}^n  \frac{k}{n-1}\ \Big[ 1 +\left( e^{(n-1) \mathrm{R}_{c} t} - 1 \right) \mathbb{E}\left(\sigma_1\left( \mathcal{A}_{i, \ell}\star_M  \mathcal{X}^{(2)}_{\ell}     \right) \right) \right. \nonumber \\
  \left. +   C_{\mbox{\tiny Cher}} \left( e^{(n-1) \mathrm{R}_{c} t} - 1 \right) \Xi( \mathcal{A}_{i, \ell}\star_M  \mathcal{X}^{(2)}_{\ell}   ) \Big]  \right\}.
\end{eqnarray} 
This Lemma is proved. $\hfill \Box$

\section{Proof for Genearlized Hanson-Wright Inequality}\label{sec:Proof for Genearlized Hanson-Wright Inequality}

We have prepared all required ingredients to prove the main theorem, Theorem~\ref{thm:GHW Inequality}, of this work. 

\textbf{Proof:}
From Eqs.~\eqref{eq:tail prob decomp1}~\eqref{eq:tail prob decomp for each j}, Lemma~\ref{lma:diag part tail bound} and Lemma~\ref{lma:coup part tail bound}, this theorem is proved.
$\hfill \Box$

\section{Conclusion}\label{sec:Conclusion} 

In this work, we extend the Hanson-Wright inequality from the quadratic sum of independent random variables to the Hanson-Wright inequality for the Ky Fan $k$-norm for the polynomial function of the quadratic sum of random tensors under Einstein product. We separate the quadratic tensors sum into the diagonal part and the coupling part. For the diagonal part, the generalized tensor Chernoff bound from~\cite{HW_T_SYChang_2021} is applied directly since each term in the diagonal part is independent of each other. For the coupling part, we apply the decoupling inequality to obtain the tail bound for the coupling part by introducing independent copies of random tensors. Then, we can apply the generalized tensor Chernoff bound again to get the tail probability of the Ky Fan $k$-norm of the coupling sum of independent random tensors. Finally, the generalized Hanson-Wright inequality for the Ky Fan $k$-norm for the polynomial function of the quadratic sum of random tensors can be obtained by the combination of the bound from the diagonal sum part and the bound from the coupling sum part.\\\\

\bibliographystyle{IEEETran}
\bibliography{Gen_HW_E_Prod_Bib}

\end{document}